# An Information-Spectrum Approach to Joint Source-Channel Coding [†]

Te Sun HAN[‡]

May 6, 2000

**Abstract:** Given a general source $\mathbf{V} = \{V^n\}_{n=1}^\infty$ with *countably infinite* source alphabet and a general channel $\mathbf{W} = \{W^n\}_{n=1}^\infty$ with arbitrary *abstract* channel input and output alphabets, we study the joint source-channel coding problem from the information-spectrum point of view. First, we generalize Feinstein's lemma (direct part) and Verdú-Han's lemma (converse part) so as to be applicable to the general joint source-channel coding problem. Based on these lemmas, we establish a sufficient condition as well as a necessary condition for the source $\mathbf{V}$ to be reliably transmissible over the channel $\mathbf{W}$ with asymptotically vanishing probability of error. It is shown that our sufficient condition coincides with the sufficient condition derived by Vembu, Verdú and Steinberg, whereas our necessary condition is much stronger than the necessary condition derived by them. Actually, our necessary condition coincide with our sufficient condition if we disregard some asymptotically vanishing terms appearing in those conditions. Also, it is shown that *Separation Theorem* in the generalized sense always holds. In addition, we demonstrate a sufficient condition as well as a necessary condition for the $\varepsilon$-transmissibility ($0 \leq \varepsilon < 1$). Finally, the separation theorem of the traditional standard form is shown to hold for the class of sources and channels that satisfy the ( semi-) strong converse property.

**Index terms:** general source, general channel, joint source-channel coding, separation theorem, information-spectrum, transmissibility, generalized Feinstein's lemma, generalized Verdú-Han's lemma

# 1 Introduction

Given a source $\mathbf{V} = \{V^n\}_{n=1}^{\infty}$ and a channel $\mathbf{W} = \{W^n\}_{n=1}^{\infty}$, the *joint source-channel coding* means that the encoder maps the output from the source directly to the channel input (*one step encoding*), where the probability of decoding error is required to vanish as block-length $n$ tends to $\infty$. In usual situations, however, the joint source-channel coding can be decomposed into separate *source coding* and *channel coding* (*two step encoding*). This two step encoding does not cause any disadvantages from the standpoint of asymptotically vanishing error probabilities, provided that the so-called *Separation Theorem* holds.

Typically, the traditional separation theorem, which we call the separation theorem in the *narrow sense*, states that if the infimum $R_f(\mathbf{V})$ of all achievable fixed-length coding rates for the source $\mathbf{V}$ is smaller than the capacity $C(\mathbf{W})$ for the channel $\mathbf{W}$ then the source $\mathbf{V}$ is reliably transmissible by two step encoding over the channel $\mathbf{W}$; whereas if $R_f(\mathbf{V})$ is larger than $C(\mathbf{W})$ then the reliable transmission is impossible. While the former statement is always true for any general source $\mathbf{V}$ and any general channel $\mathbf{W}$, the latter statement is *not* always true. Then, a very natural question may be raised for what class of sources and channels and in what sense the separation theorem holds in general.

Shannon [1] has first shown that the separation theorem holds for the class of stationary memoryless sources and channels. Since then, this theorem has received extensive attention by a number of researchers who have attempted to prove versions that apply to more and more general classes of sources and channels. Among others, for example, Dobrushin [2] and Hu [4] have studied the separation theorem problem in the context of information-stable sources and channels.

Recently, on the other hand, Vembu, Verdú and Steinberg [5] have put this problem in a much more general information-spectrum context with any general source $\mathbf{V}$ and any general channel $\mathbf{W}$. From the viewpoint of information spectra, they have generalized the notion of the separation theorem and shown that, in many cases even with $R_f(\mathbf{V}) > C(\mathbf{W})$, it is possible to reliably transmit the output of the source $\mathbf{V}$ over the channel $\mathbf{W}$. Furthermore, in terms of information spectra, they have established a sufficient condition for the transmissibility as well as a necessary condition for the transmissibility. It should be noticed here that, in general joint source-channel coding situations, what indeed matters is not the validity problem of the separation theorem but what is the necessary and sufficient condition for the transmissibility. However, while their sufficient condition looks simple and significantly tight, their necessary condition is very far from tight.

The present paper was mainly motivated by the reasonable question why these two conditions are very far from one another. In Section 3, we first demonstrate two fundamental lemmas: a generalization of Feinstein's lemma [11] and a generalization of Verdú-Han's lemma [8], which provide with the firm basis for the key results to be stated in the subsequent sections. These lemmas are of the information-spectrum forms in nice accordance with the general joint source-channel coding with *countably infinite* source alphabet, arbitrary *abstract* channel input and output alphabets. In Section 4, given a general source **V** and a general channel **W** we establish, in terms of information-spectra, a sufficient condition (*Direct theorem*) for the transmissibility as well as a necessary condition (*Converse theorem*) for the transmissibility. These two conditions are very close from each other, and actually coincides with one another if we disregard some asymptotically vanishing term. In this sense, we may regard these conditions together as specifying a "necessary and sufficient condition" for the transmissibility.

Next, we equivalently rewrite these conditions in the forms useful to see the relation with the separation theorem. As a consequence, it turns out that the equivalent form of our sufficient condition just coincides with the sufficient condition given by Vembu, Verdú and Steinberg [5], whereas the equivalent form of our necessary condition is much stronger than the necessary condition given by them. In this connection, one of our main conclusions is that *Separation Thorem* in the generalized sense always holds for all the sources and channels that satisfy this equivalent sufficient condition.

On the other hand, in Section 5, we demonstrate a sufficient condition as well as a necessary condition for the $\varepsilon$-transmissibility, which is the generalization of the sufficient condition as well as the necessary condition as was shown in Section 4. Finally, in Section 6, we restrict the class of sources and channels to those that satisfy the strong converse property (or the semi-strong converse property) to show that the separation theorem in the narrow sense holds for this class. This theorem corresponds to the standarad typical separation theorem as was described in [5, Theorem 3].

# 2 Basic Notation and Definitions

In this preliminary section, we prepare the basic notation and definitions which will be used in the subsequent sections.

## 2.1 General Source

Let us first give here the formal defintion of the general source. The *general source* is defined as an infinite sequence $\mathbf{V} = \{V^n = (V_1^{(n)}, \cdots, V_n^{(n)})\}_{n=1}^{\infty}$ of $n$-dimensional random variables $V^n$ where each component random variable $V_i^{(n)}$ $(1 \leq i \leq n)$ takes values in a *countably infinite* set $\mathcal{V}$ that we call the *source alphabet*. It should be noted here that each component of $V^n$ may change depending on block length $n$. This implies that the sequence $\mathbf{V}$ is quite general in the sense that it may not satisfy even the consistency condition as usual processes, where the consistency condition means that for any integers $m, n$ such that $m < n$ it holds that $V_i^{(m)} \equiv V_i^{(n)}$ for all $i = 1, 2, \cdots, m$. The class of sources thus defined covers a very wide range of sources including all nonstationary and/or nonergodic sources (cf. Han and Verdú [6]).

## 2.2 General Channel

The formal definition of the general channel is as follows. Let $\mathcal{X}, \mathcal{Y}$ be arbitrary *abstract* (*not* necessarily countable) sets, which we call the *input alphabet* and the *output alphabet*, respectively. The *general channel* is defined as an infinite sequence $\mathbf{W} = \{W^n(\cdot|\cdot) : \mathcal{X}^n \to \mathcal{Y}^n\}_{n=1}^{\infty}$ of $n$-dimensional probability transition matrices $W^n$, where $W^n(\mathbf{y}|\mathbf{x})$ ($\mathbf{x} \in \mathcal{X}^n, \mathbf{y} \in \mathcal{Y}^n$) denotes the conditonal probability of $\mathbf{y}$ given $\mathbf{x}$.[*] The class of channels thus defined covers a very wide range of channels including all nonstationary and/or nonergodic channels with arbitrary memory structures (cf. Han and Verdú [6]).

**Remark 2.1** A more reasonable definition of the general source is the following. Let $\{\mathcal{V}_n\}_{n=1}^{\infty}$ be any sequence of *arbitrary* source alphabets $\mathcal{V}_n$ (a countabley infinite or abstract set) and let $V_n$ be any random variable taking values in $\mathcal{V}_n$ ($n = 1, 2, \cdots$). Then, the sequence $\mathbf{V} = \{V_n\}_{n=1}^{\infty}$ of random variables $V_n$ is called a *general source* (cf. Verdú and Han [7]). The above definition is a special case of this general source with $\mathcal{V}_n = \mathcal{V}^n$ ($n = 1, 2, \cdots$).

---

[*]In the case where the output alphabet $\mathcal{Y}$ is *abstract*, $W^n(\mathbf{y}|\mathbf{x})$ is understood to be the (conditional) probability measure element $W^n(d\mathbf{y}|\mathbf{x})$ that is measurable in $\mathbf{x}$.

On the other hand, a more reasonable definition of the general channel is the following. Let $\{W_n : \mathcal{X}_n \to \mathcal{Y}_n\}_{n=1}^\infty$ be any sequence of *arbitrary* probability transition matrices, where $\mathcal{X}_n, \mathcal{Y}_n$ are arbitrary abstract sets. Then, the sequence $\mathbf{W} = \{W_n\}_{n=1}^\infty$ of probability transition matrices $W_n$ is called a *general channel* (cf. Han [9]). The above definition is a special case of this general channel with $\mathcal{X}_n = \mathcal{X}^n, \mathcal{Y}_n = \mathcal{Y}^n$ ($n = 1, 2, \cdots$).

The key results in this paper (Lemma 3.1, Lemma 3.2, Theorem 4.1, Theorem 4.2, Theorem 4.3, Theorem 4.4, Theorem 5.1, Theorem 5.2 and Theorem 6.5 ) continue to be valid as well also in this more general setting with $\mathcal{V}^n, V^n, \mathbf{V}$ and $\mathcal{X}^n, \mathcal{Y}^n, W^n, \mathbf{W}$ replaced by $\mathcal{V}_n, V_n, \mathbf{V}$ and $\mathcal{X}_n, \mathcal{Y}_n, W_n, \mathbf{W}$, respectively.

In the sequel we use the convention that $P_Z(\cdot)$ denotes the probability distribution of a random variable $Z$, whereas $P_{Z|U}(\cdot|\cdot)$ denotes the conditional probability distribution of a random variable $Z$ given a random variable $U$. □

## 2.3 Joint Source-Channel Coding

Let $\mathbf{V} = \{V^n = (V_1^{(n)}, \cdots, V_n^{(n)})\}_{n=1}^\infty$ be any general source, and let $\mathbf{W} = \{W^n(\cdot|\cdot) : \mathcal{X}^n \to \mathcal{Y}^n\}_{n=1}^\infty$ be any general channel. We consider an *encoder* $\varphi_n : \mathcal{V}^n \to \mathcal{X}^n$ and an *decoder* $\psi_n : \mathcal{Y}^n \to \mathcal{V}^n$, and put $X^n = \varphi_n(V^n)$. Then, denoting by $Y^n$ the output from the channel $W^n$ due to the input $X^n$, we have the obvious relation:

$$V^n \to X^n \to Y^n \quad \text{(a Markov chain)}. \tag{2.1}$$

The probability $\varepsilon_n$ of *decoding error* with the code $(\varphi_n, \psi_n)$ is defined by

$$\begin{aligned}\varepsilon_n &\equiv \Pr\{V^n \neq \psi_n(Y^n)\} \\ &= \sum_{\mathbf{v} \in \mathcal{V}^n} P_{V^n}(\mathbf{v}) W^n(\mathcal{D}^c(\mathbf{v})|\varphi_n(\mathbf{v})),\end{aligned} \tag{2.2}$$

where $\mathcal{D}(\mathbf{v}) \equiv \{\mathbf{y} \in \mathcal{Y}^n | \psi_n(\mathbf{y}) = \mathbf{v}\}$ ($\forall \mathbf{v} \in \mathcal{V}^n$) ($\mathcal{D}(\mathbf{v})$ is called the *decoding set* for $\mathbf{v}$) and "$c$" denotes the complement of a set. A pair $(\varphi_n, \psi_n)$ with probability $\varepsilon_n$ of decoding error is simply called a joint source-channel $(n, \varepsilon_n)$ code.

We now define the *transmissibility* in terms of joint source-channel codes $(n, \varepsilon_n)$ as

**Definition 2.1**

Source **V** is transmissible over channel **W** $\overset{\text{def}}{\iff}$ There exists an $(n, \varepsilon_n)$ code such that $\lim_{n \to \infty} \varepsilon_n = 0$.

With this definition of transmissibility, in the following sections we shall establish a sufficient condition as well as a necessary condition for the transmissibility when we are given a geneal source **V** and a general channel **W**. These two conditions are *very close* to each other and can actually be seen as giving *the same* condition if we disregard an asymptotically negligible term $\gamma_n \to 0$ appearing in those conditions (cf. Section 4).

## 3 Fundamental Lemmas

In this section, we prepare two fundamental lemmas that are needed in the next section in order to establish the main theorems (*Direct part* and *Converse part*).

**Lemma 3.1 (Generalization of Feinstein's lemma)** Given a general source $\mathbf{V} = \{V^n\}_{n=1}^{\infty}$ and a general channel $\mathbf{W} = \{W^n\}_{n=1}^{\infty}$, let $X^n$ be any input random variable taking values in $\mathcal{X}^n$, which may be arbitrarily correlated to the source variable $V^n$, and $Y^n$ be the channel output via $W^n$ due to the channel input $X^n$. Then, for every $n = 1, 2, \cdots$, there exists an $(n, \varepsilon_n)$ code such that

$$\varepsilon_n \leq \Pr\left\{ \frac{1}{n} \log \frac{W^n(Y^n|X^n)}{P_{Y^n}(Y^n)} \leq \frac{1}{n} \log \frac{1}{P_{V^n}(V^n)} + \gamma \right\} + e^{-n\gamma}, \quad (3.1)$$

where[†] $\gamma > 0$ is an arbitrary positive number.

**Remark 3.1** In a special case where the source $\mathbf{V} = \{V^n\}_{n=1}^{\infty}$ is uniformly distributed on the massage set $\mathcal{M}_n = \{1, 2, \cdots, M_n\}$, it follows that

$$\frac{1}{n} \log \frac{1}{P_{V^n}(V^n)} = \frac{1}{n} \log M_n,$$

---

[†]In the case where the input and output alphabets $\mathcal{X}, \mathcal{Y}$ are *abstract* (*not* necessarily countable), $\frac{W^n(Y^n|X^n)}{P_{Y^n}(Y^n)}$ in (3.1) is understood to be $g(Y^n|X^n)$, where $g(\mathbf{y}|\mathbf{x}) \equiv \frac{W^n(d\mathbf{y}|\mathbf{x})}{P_{Y^n}(d\mathbf{y})}$ $= \frac{W^n(d\mathbf{y}|\mathbf{x})P_{X^n}(d\mathbf{x})}{P_{Y^n}(d\mathbf{y})P_{X^n}(d\mathbf{x})} = \frac{P_{X^nY^n}(d\mathbf{x},d\mathbf{y})}{P_{X^n}(d\mathbf{x})P_{Y^n}(d\mathbf{y})}$ is the Radon-Nikodym derivative that is measurable in $(\mathbf{x}, \mathbf{y})$.

which implies that the entropy spectrum[‡] of the source $\mathbf{V} = \{V^n\}_{n=1}^{\infty}$ is exactly one point spectrum concentrated on $\frac{1}{n}\log M_n$. Therefore, in this special case, Lemma 3.1 reduecs to Feinstein's lemma [11]. □

*Proof of Lemma3.1*:

For each $\mathbf{v} \in \mathcal{V}^n$, generate $\mathbf{x}(\mathbf{v}) \in \mathcal{X}^n$ at random according to the conditional distribution $P_{X^n|V^n}(\cdot|\mathbf{v})$ and let $\mathbf{x}(\mathbf{v})$ be the codeword for $\mathbf{v}$. In other words, we define the encoder $\varphi_n : \mathcal{V}^n \to \mathcal{X}^n$ as $\mathbf{x}(\mathbf{v}) = \varphi_n(\mathbf{v})$, where $\{\mathbf{x}(\mathbf{v}) \mid \forall \mathbf{v} \in \mathcal{V}^n\}$ are all independently generated. We define the decoder $\psi_n : \mathcal{Y}^n \to \mathcal{V}^n$ as follows: Set

$$S_n = \left\{(\mathbf{v},\mathbf{x},\mathbf{y}) \in \mathcal{Z}^n \left| \frac{1}{n}\log\frac{W^n(\mathbf{y}|\mathbf{x})}{P_{Y^n}(\mathbf{y})} > \frac{1}{n}\log\frac{1}{P_{V^n}(\mathbf{v})} + \gamma \right.\right\}, \quad (3.2)$$

$$S_n(\mathbf{v}) = \{(\mathbf{x},\mathbf{y}) \in \mathcal{X}^n \times \mathcal{Y}^n \mid (\mathbf{v},\mathbf{x},\mathbf{y}) \in S_n\}, \quad (3.3)$$

where for simplicity we have put $\mathcal{Z}^n \equiv \mathcal{V}^n \times \mathcal{X}^n \times \mathcal{Y}^n$. Suppose that the decoder $\psi_n$ received a channel output $\mathbf{y} \in \mathcal{Y}^n$. If there exists one and only one $\mathbf{v} \in \mathcal{V}^n$ such that $(\mathbf{x}(\mathbf{v}), \mathbf{y}) \in S_n(\mathbf{v})$, define the decoder as $\mathbf{v} = \psi_n(\mathbf{y})$ with that $\mathbf{v}$; otherwise, let the output of the decoder $\psi_n(\mathbf{y}) \in \mathcal{V}^n$ be arbitrary. Then, the probability $\overline{\varepsilon}_n$ of error for this pair $(\varphi_n, \psi_n)$ of encoder and decoder (averaged over all the realizatioins of the random code) is given by

$$\overline{\varepsilon}_n = \sum_{\mathbf{v} \in \mathcal{V}^n} P_{V^n}(\mathbf{v})\overline{\varepsilon}_n(\mathbf{v}), \quad (3.4)$$

where $\overline{\varepsilon}_n(\mathbf{v})$ is the probability of error (averaged over all the realizatioins of the random code) when $\mathbf{v} \in \mathcal{V}^n$ is the source output. We can evaluate $\overline{\varepsilon}_n(\mathbf{v})$ as

$$\overline{\varepsilon}_n(\mathbf{v}) \leq \Pr\{(\mathbf{x}(\mathbf{v}), Y^n) \notin S_n(\mathbf{v})\}$$
$$+ \Pr\left\{\bigcup_{\mathbf{v}':\mathbf{v}'\neq\mathbf{v}}\{(\mathbf{x}(\mathbf{v}'), Y^n) \in S_n(\mathbf{v}')\}\right\}$$

---

[‡]The probablity distribution of $\frac{1}{n}\log\frac{1}{P_{V^n}(V^n)}$ is called the *entropy spectrum* of the source $\mathbf{V} = \{V^n\}_{n=1}^{\infty}$, whereas the probability distribution of $\frac{1}{n}\log\frac{W^n(Y^n|X^n)}{P_{Y^n}(Y^n)}$ is called the *mutual information spectrum* of the channel $\mathbf{W} = \{W^n\}_{n=1}^{\infty}$ given the input $\mathbf{X} = \{X^n\}_{n=1}^{\infty}$ (cf. Han and Verdú [6]).

$$\begin{aligned}
&\leq \Pr\{(\mathbf{x}(\mathbf{v}), Y^n) \notin S_n(\mathbf{v})\} \\
&\quad + \sum_{\mathbf{v}':\mathbf{v}'\neq\mathbf{v}} \Pr\{(\mathbf{x}(\mathbf{v}'), Y^n) \in S_n(\mathbf{v}')\},
\end{aligned} \quad (3.5)$$

where $Y^n$ is the channel output via $W^n$ due to the channel input $\mathbf{x}(\mathbf{v})$. The first term on the right-hand side of (3.5) is written as

$$\begin{aligned}
A_n(\mathbf{v}) &\equiv \Pr\{(\mathbf{x}(\mathbf{v}), Y^n) \notin S_n(\mathbf{v})\} \\
&= \sum_{(\mathbf{x},\mathbf{y})\notin S_n(\mathbf{v})} P_{X^n Y^n | V^n}(\mathbf{x}, \mathbf{y} | \mathbf{v}).
\end{aligned}$$

Hence,

$$\begin{aligned}
\sum_{\mathbf{v}\in\mathcal{V}^n} P_{V^n}(\mathbf{v}) A_n(\mathbf{v}) &= \sum_{\mathbf{v}\in\mathcal{V}^n} P_{V^n}(\mathbf{v}) \sum_{(\mathbf{x},\mathbf{y})\notin S_n(\mathbf{v})} P_{X^n Y^n | V^n}(\mathbf{x}, \mathbf{y} | \mathbf{v}) \\
&= \sum_{(\mathbf{v},\mathbf{x},\mathbf{y})\notin S_n} P_{V^n X^n Y^n}(\mathbf{v}, \mathbf{x}, \mathbf{y}) \\
&= \Pr\{V^n X^n Y^n \notin S_n\}.
\end{aligned} \quad (3.6)$$

On the other hand, noting that $\mathbf{x}(\mathbf{v}'), \mathbf{x}(\mathbf{v})$ ($\mathbf{v}' \neq \mathbf{v}$) are independent and hence $\mathbf{x}(\mathbf{v}'), Y^n$ are also independent, the second term on the right-hand side of (3.5) is evaluated as

$$\begin{aligned}
B_n(\mathbf{v}) &\equiv \sum_{\mathbf{v}':\mathbf{v}'\neq\mathbf{v}} \Pr\{(\mathbf{x}(\mathbf{v}'), Y^n) \in S_n(\mathbf{v}')\} \\
&= \sum_{\mathbf{v}':\mathbf{v}'\neq\mathbf{v}} \sum_{(\mathbf{x},\mathbf{y})\in S_n(\mathbf{v}')} P_{Y^n | V^n}(\mathbf{y} | \mathbf{v}) P_{X^n | V^n}(\mathbf{x} | \mathbf{v}') \\
&\leq \sum_{\mathbf{v}'\in\mathcal{V}^n} \sum_{(\mathbf{x},\mathbf{y})\in S_n(\mathbf{v}')} P_{Y^n | V^n}(\mathbf{y} | \mathbf{v}) P_{X^n | V^n}(\mathbf{x} | \mathbf{v}').
\end{aligned}$$

Hence,

$$\begin{aligned}
&\sum_{\mathbf{v}\in\mathcal{V}^n} P_{V^n}(\mathbf{v}) B_n(\mathbf{v}) \\
&\leq \sum_{\mathbf{v}\in\mathcal{V}^n} \sum_{\mathbf{v}'\in\mathcal{V}^n} \sum_{(\mathbf{x},\mathbf{y})\in S_n(\mathbf{v}')} P_{V^n}(\mathbf{v}) P_{Y^n | V^n}(\mathbf{y} | \mathbf{v}) P_{X^n | V^n}(\mathbf{x} | \mathbf{v}') \\
&= \sum_{\mathbf{v}'\in\mathcal{V}^n} \sum_{(\mathbf{x},\mathbf{y})\in S_n(\mathbf{v}')} P_{Y^n}(\mathbf{y}) P_{X^n | V^n}(\mathbf{x} | \mathbf{v}').
\end{aligned} \quad (3.7)$$

Since, if $(\mathbf{x}, \mathbf{y}) \in S_n(\mathbf{v}')$, then by means of (3.2), (3.3) we have

$$P_{Y^n}(\mathbf{y}) \leq P_{V^n}(\mathbf{v}') W^n(\mathbf{y}|\mathbf{x}) e^{-n\gamma},$$

(3.7) is further transformed to

$$\sum_{\mathbf{v} \in \mathcal{V}^n} P_{V^n}(\mathbf{v}) B_n(\mathbf{v})$$
$$\leq e^{-n\gamma} \sum_{\mathbf{v}' \in \mathcal{V}^n} \sum_{(\mathbf{x},\mathbf{y}) \in S_n(\mathbf{v}')} P_{V^n}(\mathbf{v}') P_{X^n|V^n}(\mathbf{x}|\mathbf{v}') W^n(\mathbf{y}|\mathbf{x})$$
$$\leq e^{-n\gamma} \sum_{(\mathbf{v}',\mathbf{x},\mathbf{y}) \in \mathcal{Z}^n} P_{V^n}(\mathbf{v}') P_{X^n|V^n}(\mathbf{x}|\mathbf{v}') W^n(\mathbf{y}|\mathbf{x})$$
$$= e^{-n\gamma}. \qquad (3.8)$$

Then, from (3.4), (3.6) and (3.8) it follows that

$$\begin{aligned}\overline{\varepsilon}_n &= \sum_{\mathbf{v} \in \mathcal{V}^n} P_{V^n}(\mathbf{v}) \overline{\varepsilon}_n(\mathbf{v}) \\ &\leq \sum_{\mathbf{v} \in \mathcal{V}^n} P_{V^n}(\mathbf{v}) A_n(\mathbf{v}) + \sum_{\mathbf{v} \in \mathcal{V}^n} P_{V^n}(\mathbf{v}) B_n(\mathbf{v}) \\ &\leq \Pr\{V^n X^n Y^n \notin S_n\} + e^{-n\gamma}.\end{aligned}$$

Thus, there must exist a deterministic $(n, \varepsilon_n)$ code such that

$$\varepsilon_n \leq \Pr\{V^n X^n Y^n \notin S_n\} + e^{-n\gamma},$$

thereby proving Lemma 3.1. $\square$

**Lemma 3.2 (Generalization of Verdú-Han's lemma)** Let $\mathbf{V} = \{V^n\}_{n=1}^{\infty}$ and $\mathbf{W} = \{W^n\}_{n=1}^{\infty}$ be a general source and a general channel, respectively, and let $\varphi_n : \mathcal{V}^n \to \mathcal{X}^n$ be the encoder of an $(n, \varepsilon_n)$ code for $(V^n, W^n)$. Put $X^n = \varphi_n(V^n)$ and let $Y^n$ be the channel output via $W^n$ due to the channel input $X^n$. Then, for every $n = 1, 2, \cdots$, it holds that

$$\varepsilon_n \geq \Pr\left\{\frac{1}{n} \log \frac{W^n(Y^n|X^n)}{P_{Y^n}(Y^n)} \leq \frac{1}{n} \log \frac{1}{P_{V^n}(V^n)} - \gamma\right\} - e^{-n\gamma}, \qquad (3.9)$$

where $\gamma > 0$ is an arbitrary positive number.

**Remark 3.2** In a special case where the source $\mathbf{V} = \{V^n\}_{n=1}^{\infty}$ is uniformly distributed on the massage set $\mathcal{M}_n = \{1, 2, \cdots, M_n\}$, it follows that

$$\frac{1}{n} \log \frac{1}{P_{V^n}(V^n)} = \frac{1}{n} \log M_n,$$

which implies that the entropy spectrum of the source $\mathbf{V} = \{V^n\}_{n=1}^{\infty}$ is exactly one point spectrum concentrated on $\frac{1}{n}\log M_n$. Therefore, in this special case, Lemma 3.2 reduecs to Verdú-Han's lemma [8]. □

*Proof of Lemma3.2*

Define
$$L_n = \left\{(\mathbf{v}, \mathbf{x}, \mathbf{y}) \in \mathcal{Z}^n \left| \frac{1}{n}\log\frac{W^n(\mathbf{y}|\mathbf{x})}{P_{Y^n}(\mathbf{y})} \leq \frac{1}{n}\log\frac{1}{P_{V^n}(\mathbf{v})} - \gamma \right.\right\}, \quad (3.10)$$
and, for each $\mathbf{v} \in \mathcal{V}^n$, set
$$\mathcal{D}(\mathbf{v}) = \{\mathbf{y} \in \mathcal{Y}^n | \psi_n(\mathbf{y}) = \mathbf{v}\},$$
that is, $\mathcal{D}(\mathbf{v})$ is the decoding set for $\mathbf{v}$. Moreover, for each $(\mathbf{v}, \mathbf{x}) \in \mathcal{V}^n \times \mathcal{X}^n$, set
$$\mathcal{B}(\mathbf{v}, \mathbf{x}) = \{\mathbf{y} \in \mathcal{Y}^n | (\mathbf{v}, \mathbf{x}, \mathbf{y}) \in L_n\}. \quad (3.11)$$
Then, noting the Markov chain property (2.1), we have
$$\begin{aligned}
&\Pr\{V^n X^n Y^n \in L_n\} \\
&= \sum_{(\mathbf{v},\mathbf{x},\mathbf{y}) \in L_n} P_{V^n X^n Y^n}(\mathbf{v}, \mathbf{x}, \mathbf{y}) \\
&= \sum_{(\mathbf{v},\mathbf{x}) \in \mathcal{V}^n \times \mathcal{X}^n} P_{V^n X^n}(\mathbf{v}, \mathbf{x}) W^n(\mathcal{B}(\mathbf{v}, \mathbf{x})|\mathbf{x}) \\
&= \sum_{(\mathbf{v},\mathbf{x}) \in \mathcal{V}^n \times \mathcal{X}^n} P_{V^n X^n}(\mathbf{v}, \mathbf{x}) W^n(\mathcal{B}(\mathbf{v}, \mathbf{x}) \cap \mathcal{D}^c(\mathbf{v})|\mathbf{x}) \\
&\quad + \sum_{(\mathbf{v},\mathbf{x}) \in \mathcal{V}^n \times \mathcal{X}^n} P_{V^n X^n}(\mathbf{v}, \mathbf{x}) W^n(\mathcal{B}(\mathbf{v}, \mathbf{x}) \cap \mathcal{D}(\mathbf{v})|\mathbf{x}) \\
&\leq \sum_{(\mathbf{v},\mathbf{x}) \in \mathcal{V}^n \times \mathcal{X}^n} P_{V^n X^n}(\mathbf{v}, \mathbf{x}) W^n(\mathcal{D}^c(\mathbf{v})|\mathbf{x}) \\
&\quad + \sum_{(\mathbf{v},\mathbf{x}) \in \mathcal{V}^n \times \mathcal{X}^n} P_{V^n X^n}(\mathbf{v}, \mathbf{x}) W^n(\mathcal{B}(\mathbf{v}, \mathbf{x}) \cap \mathcal{D}(\mathbf{v})|\mathbf{x}) \\
&= \varepsilon_n + \sum_{(\mathbf{v},\mathbf{x}) \in \mathcal{V}^n \times \mathcal{X}^n} P_{V^n X^n}(\mathbf{v}, \mathbf{x}) W^n(\mathcal{B}(\mathbf{v}, \mathbf{x}) \cap \mathcal{D}(\mathbf{v})|\mathbf{x}) \\
&= \varepsilon_n + \sum_{(\mathbf{v},\mathbf{x}) \in \mathcal{V}^n \times \mathcal{X}^n} P_{V^n X^n}(\mathbf{v}, \mathbf{x}) \sum_{\mathbf{y} \in \mathcal{B}(\mathbf{v},\mathbf{x}) \cap \mathcal{D}(\mathbf{v})} W^n(\mathbf{y}|\mathbf{x}), \quad (3.12)
\end{aligned}$$

where we have used the relation:
$$\varepsilon_n = \sum_{(\mathbf{v},\mathbf{x})\in\mathcal{V}^n\times\mathcal{X}^n} P_{V^n X^n}(\mathbf{v},\mathbf{x}) W^n(\mathcal{D}^c(\mathbf{v})|\mathbf{x}).$$

Now, it follows from (3.10) and (3.11) that $\mathbf{y} \in \mathcal{B}(\mathbf{v},\mathbf{x})$ implies
$$W^n(\mathbf{y}|\mathbf{x}) \le \frac{e^{-n\gamma} P_{Y^n}(\mathbf{y})}{P_{V^n}(\mathbf{v})},$$

which is substituted into the right-hand side of (3.12) to yield
$$\begin{aligned}
&\Pr\{V^n X^n Y^n \in L_n\} \\
&\le \varepsilon_n + e^{-n\gamma} \sum_{(\mathbf{v},\mathbf{x})\in\mathcal{V}^n\times\mathcal{X}^n} P_{X^n|V^n}(\mathbf{x}|\mathbf{v}) \sum_{\mathbf{y}\in\mathcal{B}(\mathbf{v},\mathbf{x})\cap\mathcal{D}(\mathbf{v})} P_{Y^n}(\mathbf{y}) \\
&\le \varepsilon_n + e^{-n\gamma} \sum_{(\mathbf{v},\mathbf{x})\in\mathcal{V}^n\times\mathcal{X}^n} P_{X^n|V^n}(\mathbf{x}|\mathbf{v}) P_{Y^n}(\mathcal{D}(\mathbf{v})) \\
&= \varepsilon_n + e^{-n\gamma} \sum_{\mathbf{v}\in\mathcal{V}^n} P_{Y^n}(\mathcal{D}(\mathbf{v})) \\
&= \varepsilon_n + e^{-n\gamma},
\end{aligned}$$

thereby proving the claim of the lemma. □

## 4 Theorems on Transmissibility

In this section we give both of a sufficient condition and a necessary condition for the transmissibility with a given general souce $\mathbf{V} = \{V^n\}_{n=1}^\infty$ and a given general channel $\mathbf{W} = \{W^n\}_{n=1}^\infty$.

First, Lemma 3.1 immediately leads us to the following direct theorem:

**Theorem 4.1 (Direct theorem)** Let $\mathbf{V} = \{V^n\}_{n=1}^\infty$, $\mathbf{W} = \{W^n\}_{n=1}^\infty$ be a general source and a general channel, respectively. If there exist some channel input $\mathbf{X} = \{X^n\}_{n=1}^\infty$, which may be arbitrarily correlated to the output of the source $\mathbf{V} = \{V^n\}_{n=1}^\infty$, and also some sequence $\{\gamma_n\}_{n=1}^\infty$ of real nubers with
$$\gamma_n > 0,\ \gamma_n \to 0 \text{ and } n\gamma_n \to \infty \quad (n \to \infty) \tag{4.1}$$
for which it holds that
$$\lim_{n\to\infty} \Pr\left\{\frac{1}{n}\log\frac{W^n(Y^n|X^n)}{P_{Y^n}(Y^n)} \le \frac{1}{n}\log\frac{1}{P_{V^n}(V^n)} + \gamma_n\right\} = 0, \tag{4.2}$$

then the source $\mathbf{V} = \{V^n\}_{n=1}^{\infty}$ is transmissible over the channel $\mathbf{W} = \{W^n\}_{n=1}^{\infty}$, where $Y^n$ is the channel output via $W^n$ due to the channel input $X^n$.

*Proof:*

Since in Lemma 3.1 we can choose the constant $\gamma > 0$ so as to depend on $n$, let us take, instead of $\gamma$, an arbitrary $\{\gamma_n\}_{n=1}^{\infty}$ satisfying condition (4.1). Then, the second term on the right-hand side of (3.1) vanishes as $n$ tends to $\infty$, and hence it follows from (4.2) that the right-hand side of (3.1) vanishes as $n$ tends to $\infty$. Therefore, the $(n, \varepsilon_n)$ code as specified in Lemma 3.1 satisfies $\lim_{n\to\infty} \varepsilon_n = 0$. □

Next, Lemma 3.2 immediately leads us to the following converse theorem:

**Theorem 4.2 (Converse theorem)** Suppose that a general source $\mathbf{V} = \{V^n\}_{n=1}^{\infty}$ is transmissible over a general channel $\mathbf{W} = \{W^n\}_{n=1}^{\infty}$. Then, for some channel input $\mathbf{X} = \{X^n\}_{n=1}^{\infty}$, which may be arbitrarily correlated to the output of the source $\mathbf{V} = \{V^n\}_{n=1}^{\infty}$, and for any sequence $\{\gamma_n\}_{n=1}^{\infty}$ satisfying condition (4.1), it holds that

$$\lim_{n\to\infty} \Pr\left\{ \frac{1}{n} \log \frac{W^n(Y^n|X^n)}{P_{Y^n}(Y^n)} \le \frac{1}{n} \log \frac{1}{P_{V^n}(V^n)} - \gamma_n \right\} = 0, \qquad (4.3)$$

where $Y^n$ is the channel output via $W^n$ due to the channel input $X^n$.

*Proof:*

If $\mathbf{V}$ is transmissible over $\mathbf{W}$, then, by Definition 2.1 there exists an $(n, \varepsilon_n)$ code such that $\lim_{n\to\infty} \varepsilon_n = 0$. Denote by $\varphi_n$ the encoder of this code and put $X^n = \varphi_n(V^n)$. Moreover, if we denote by $Y^n$ the channel output via $W^n$ due to the channel input $X^n$, then the claim of the theorem immediately follows from (3.9) in Lemma 3.2 with $\gamma_n$ instead of $\gamma$. □

**Remark 4.1** Comparing (4.3) in Theorem 4.2 with (4.2) in Theorem 4.1, we observe that the only diffence is that the sign of $\gamma_n$ is changed from $+$ to $-$. Since $\gamma_n$ vanishes as $n$ tends to $\infty$, this difference is asymptotically negligible. Thus, except for this asymptotically negligible difference, Theorem 4.1 together with Theorem 4.2 can be regarded as providing with a " necessary and sufficient condition" for the source $\mathbf{V} = \{V^n\}_{n=1}^{\infty}$ to be transmissible over the channel $\mathbf{W} = \{W^n\}_{n=1}^{\infty}$. □

Now, let us think of the implication of conditions (4.2), (4.3). First, let us think of (4.2). Putting

$$A_n = \frac{1}{n} \log \frac{W^n(Y^n|X^n)}{P_{Y^n}(Y^n)}, \quad B_n = \frac{1}{n} \log \frac{1}{P_{V^n}(V^n)}$$

for simplicity, (4.2) is written as

$$\alpha_n \equiv \Pr\{A_n \leq B_n + \gamma_n\} \to 0 \quad (n \to \infty), \tag{4.4}$$

which can be transformed to

$$\begin{aligned}
\Pr\{A_n &\leq B_n + \gamma_n\} \\
&= \sum_u \Pr\{B_n = u\} \Pr\{A_n \leq B_n + \gamma_n | B_n = u\} \\
&= \sum_u \Pr\{B_n = u\} \Pr\{A_n \leq u + \gamma_n | B_n = u\}.
\end{aligned}$$

Set

$$T_n = \{u \mid \Pr\{A_n \leq u + \gamma_n | B_n = u\} \leq \sqrt{\alpha_n}\}, \tag{4.5}$$

then by virtue of (4.4) and Markov inequality, we have

$$\Pr\{B_n \in T_n\} \geq 1 - \sqrt{\alpha_n}. \tag{4.6}$$

Let us now define the upper cumulative probabilities for $A_n, B_n$ by

$$P_n(t) = \Pr\{A_n \geq t\}, \quad Q_n(t) = \Pr\{B_n \geq t\},$$

then it follows that

$$\begin{aligned}
P_n(t) &= \sum_u \Pr\{B_n = u\} \Pr\{A_n \geq t | B_n = u\} \\
&\geq \sum_{\substack{u \in T_n: \\ u \geq t - \gamma_n}} \Pr\{B_n = u\} \Pr\{A_n \geq t | B_n = u\} \\
&\geq \sum_{\substack{u \in T_n: \\ u \geq t - \gamma_n}} \Pr\{B_n = u\} \Pr\{A_n \geq u + \gamma_n | B_n = u\}. \tag{4.7}
\end{aligned}$$

On the other hand, by means of (4.5), $u \in T_n$ implies that

$$\Pr\{A_n \geq u + \gamma_n | B_n = u\} \geq 1 - \sqrt{\alpha_n}.$$

Theore, by (4.6), (4.7) it is concluded that

$$\begin{aligned} P_n(t) &\geq (1 - \sqrt{\alpha_n}) \sum_{\substack{u \in T_n: \\ u \geq t - \gamma_n}} \Pr\{B_n = u\} \\ &\geq (1 - \sqrt{\alpha_n})(Q_n(t - \gamma_n) - \Pr\{B_n \notin T_n\}) \\ &\geq (1 - \sqrt{\alpha_n})(Q_n(t - \gamma_n) - \sqrt{\alpha_n}) \\ &\geq Q_n(t - \gamma_n) - 2\sqrt{\alpha_n}. \end{aligned}$$

This means that, for all $t$, the upper cumulative probability $P_n(t)$ of $A_n$ is larger than or equal to the upper cumulative probability $Q_n(t - \gamma_n)$ of $B_n$, except for the asymptotically vanishing difference $2\sqrt{\alpha_n}$. This in turn implies that, as a whole, the mutual information spectrum of the channel is shifted to the right in comparison with the entropy spectrum of the source. With $-\gamma_n$ instead of $\gamma_n$, the same implication follows also from (4.3). This is the information-spectrum implication of the "necessary and sufficient condition" (4.2), (4.3). It is such an allocation relation between the mutual information spectrum and the entropy spectrum that enables us to make an transmissible joint source-channel coding.

However, it is not easy in general to check whether conditions (4.2), (4.3) in these forms are satisfied or not. Therefore, we consider to equivalently rewrite conditions (4.2), (4.3) into the forms easier to check. This can actually be done by re-choosing the input and output variables $X^n, Y^n$ as follows. These forms are useful in order to see the relation of conditions (4.2), (4.3) with the so-called *Separation Theorem*.

First, we show the equivalent information-spectrum form of the sufficient condition (4.2) in Theorem 4.1.

**Theorem 4.3 (Equivalence of sufficient conditions)** The following two conditions are equivalent:

1) For some channel input $\mathbf{X} = \{X^n\}_{n=1}^{\infty}$, which may be arbitrarily correlated to the output of the source $\mathbf{V} = \{V^n\}_{n=1}^{\infty}$, and for some sequence $\{\gamma_n\}_{n=1}^{\infty}$ satisfying condition (4.1), it holds that

$$\lim_{n \to \infty} \Pr\left\{ \frac{1}{n} \log \frac{W^n(Y^n|X^n)}{P_{Y^n}(Y^n)} \leq \frac{1}{n} \log \frac{1}{P_{V^n}(V^n)} + \gamma_n \right\} = 0, \qquad (4.8)$$

where $Y^n$ is the channel output via $W^n$ due to the channel input $X^n$.

2) (**Strict domination:** Vembu, Verdú and Steinberg [5]) For some channel input $\mathbf{X} = \{X^n\}_{n=1}^{\infty}$, which may be arbitrarily correlated to the output of the source $\mathbf{V} = \{V^n\}_{n=1}^\infty$, and for some sequence $\{c_n\}_{n=1}^\infty$ and some sequence $\{\gamma_n\}_{n=1}^\infty$ satisfying condition (4.1), it holds that

$$\lim_{n \to \infty} \left( \Pr\left\{ \frac{1}{n} \log \frac{1}{P_{V^n}(V^n)} \geq c_n \right\} \right.$$
$$\left. + \Pr\left\{ \frac{1}{n} \log \frac{W^n(Y^n|X^n)}{P_{Y^n}(Y^n)} \leq c_n + \gamma_n \right\} \right) = 0, \qquad (4.9)$$

where $Y^n$ is the channel output via $W^n$ due to the channel input $X^n$.

**Remark 4.2 (Separation)** The sufficient condition 2) in Theorem 4.3 means that the entropy spectrum of tha source and the mutual information spectrum of the channel are asymptotically completely split with the vacant boundary of asymptotically vanishing width $\gamma_n$, and the former is placed to the left of the latter, where these two spectra may vibrate "synchronously" with $n$. In the case where such a separation condition 2) is satisfied, we can make the transmissible joint source-channel coding in two steps as follows (*Separation* of source coding and channel coding): We first encode the source output $V^n$ at the fixed-length coding rate $c_n = \frac{1}{n} \log M_n$, and then encode the output of the source encoder into the channel. The error probabilty $\varepsilon_n$ for this two step coding is upper bounded by the sum of the average error probability of the fixed-length source coding (cf. Han [9], [10]):

$$\Pr\left\{ \frac{1}{n} \log \frac{1}{P_{V^n}(V^n)} \geq c_n \right\}$$

and the maximum error probability of the channel coding (cf. Feinstein [11]):

$$\Pr\left\{ \frac{1}{n} \log \frac{W^n(Y^n|X^n)}{P_{Y^n}(Y^n)} \leq c_n + \gamma_n \right\} + e^{-n\gamma_n}.$$

It then follows from (4.9) that both of these two error probabilities vanish as $n$ tends to $\infty$, where it should be noted that $e^{-n\gamma_n} \to 0$ as $n \to \infty$. Thus, we have $\lim_{n \to \infty} \varepsilon_n = 0$ to conclude that the source $\mathbf{V} = \{V^n\}_{n=1}^\infty$ is transmissible over the channel $\mathbf{W} = \{W^n\}_{n=1}^\infty$. This gives also another proof of Theorem 4.1. □

*Proof of Theorem 4.3:*

2) ⇒ 1): For any joint probability distribution $P_{V^n X^n}$ for $V^n$ and $X^n$, we have

$$\Pr\left\{\frac{1}{n}\log\frac{W^n(Y^n|X^n)}{P_{Y^n}(Y^n)} \le \frac{1}{n}\log\frac{1}{P_{V^n}(V^n)} + \gamma_n\right\}$$
$$\le \Pr\left\{\frac{1}{n}\log\frac{1}{P_{V^n}(V^n)} \ge c_n\right\}$$
$$+ \Pr\left\{\frac{1}{n}\log\frac{W^n(Y^n|X^n)}{P_{Y^n}(Y^n)} \le c_n + \gamma_n\right\},$$

which together with (4.9) implies (4.8).

1) ⇒ 2)   Supposing that condition 1) holds, put

$$\alpha_n \equiv \Pr\left\{\frac{1}{n}\log\frac{W^n(Y^n|X^n)}{P_{Y^n}(Y^n)} \le \frac{1}{n}\log\frac{1}{P_{V^n}(V^n)} + \gamma_n\right\}, \tag{4.10}$$

and moreover, with $\gamma'_n = \dfrac{\gamma_n}{4}$, $\delta_n = \max(\sqrt{\alpha_n}, e^{-n\gamma'_n})$, put

$$d_n = \sup\left\{R \,\Big|\, \Pr\left\{\frac{1}{n}\log\frac{1}{P_{V^n}(V^n)} \ge R\right\} > \delta_n\right\} - \gamma'_n. \tag{4.11}$$

Furthermore, if we put

$$S_n = \left\{\mathbf{v} \in \mathcal{V}^n \,\Big|\, \frac{1}{n}\log\frac{1}{P_{V^n}(\mathbf{v})} \ge d_n\right\}, \tag{4.12}$$
$$\lambda_n^{(1)} = \Pr\{V^n \in S_n\}, \quad \lambda_n^{(2)} = \Pr\{V^n \notin S_n\}, \tag{4.13}$$

then the joint probability distribution $P_{V^n X^n Y^n}$ can be written as a mixture:

$$P_{V^n X^n Y^n}(\mathbf{v}, \mathbf{x}, \mathbf{y})$$
$$= \lambda_n^{(1)} P_{\tilde{V}^n \tilde{X}^n \tilde{Y}^n}(\mathbf{v}, \mathbf{x}, \mathbf{y}) + \lambda_n^{(2)} P_{\overline{V}^n \overline{X}^n \overline{Y}^n}(\mathbf{v}, \mathbf{x}, \mathbf{y}), \tag{4.14}$$

where $P_{\tilde{V}^n \tilde{X}^n \tilde{Y}^n}$, $P_{\overline{V}^n \overline{X}^n \overline{Y}^n}$ are the conditional probability distributions of $V^n X^n Y^n$ conditioned on $V^n \in S_n$, $V^n \notin S_n$, respectively. We notice here that the Markov chain property $V^n \to X^n \to Y^n$ implies $P_{\tilde{Y}^n|\tilde{X}^n} = P_{\overline{Y}^n|\overline{X}^n}$ $= W^n$ as well as the Markov chain properties

$$\tilde{V}^n \to \tilde{X}^n \to \tilde{Y}^n, \quad \overline{V}^n \to \overline{X}^n \to \overline{Y}^n.$$

We now rewrite (4.10) as

$$\begin{aligned}\alpha_n &= \lambda_n^{(1)} \Pr\left\{\frac{1}{n}\log\frac{W^n(\tilde{Y}^n|\tilde{X}^n)}{P_{Y^n}(\tilde{Y}^n)} \leq \frac{1}{n}\log\frac{1}{P_{V^n}(\tilde{V}^n)} + \gamma_n\right\} \\ &\quad + \lambda_n^{(2)} \Pr\left\{\frac{1}{n}\log\frac{W^n(\overline{Y}^n|\overline{X}^n)}{P_{Y^n}(\overline{Y}^n)} \leq \frac{1}{n}\log\frac{1}{P_{V^n}(\overline{V}^n)} + \gamma_n\right\}.\end{aligned} \quad (4.15)$$

On the other hand, since (4.11), (4.12) lead to $\lambda_n^{(1)} > \delta_n \geq \sqrt{\alpha_n}$, it follows from (4.15) that

$$\Pr\left\{\frac{1}{n}\log\frac{W^n(\tilde{Y}^n|\tilde{X}^n)}{P_{Y^n}(\tilde{Y}^n)} \leq \frac{1}{n}\log\frac{1}{P_{V^n}(\tilde{V}^n)} + \gamma_n\right\} \leq \sqrt{\alpha_n}. \quad (4.16)$$

Then, by the definition of $\tilde{V}^n$,

$$\frac{1}{n}\log\frac{1}{P_{V^n}(\tilde{V}^n)} \geq d_n,$$

and so from (4.16), we obtain

$$\Pr\left\{\frac{1}{n}\log\frac{W^n(\tilde{Y}^n|\tilde{X}^n)}{P_{Y^n}(\tilde{Y}^n)} \leq d_n + \gamma_n\right\} \leq \sqrt{\alpha_n}. \quad (4.17)$$

Next, since it follows from (4.14) that

$$\begin{aligned}P_{Y^n}(\mathbf{y}) &= \lambda_n^{(1)} P_{\tilde{Y}^n}(\mathbf{y}) + \lambda_n^{(2)} P_{\overline{Y}^n}(\mathbf{y}) \\ &\geq \lambda_n^{(1)} P_{\tilde{Y}^n}(\mathbf{y}) \\ &\geq \delta_n P_{\tilde{Y}^n}(\mathbf{y}) \\ &\geq e^{-n\gamma_n'} P_{\tilde{Y}^n}(\mathbf{y}),\end{aligned}$$

we have

$$\frac{1}{n}\log\frac{1}{P_{Y^n}(\tilde{Y}^n)} \leq \frac{1}{n}\log\frac{1}{P_{\tilde{Y}^n}(\tilde{Y}^n)} + \gamma_n',$$

which is substituted into (4.17) to get

$$\Pr\left\{\frac{1}{n}\log\frac{W^n(\tilde{Y}^n|\tilde{X}^n)}{P_{\tilde{Y}^n}(\tilde{Y}^n)} \leq d_n + \gamma_n - \gamma_n'\right\} \leq \sqrt{\alpha_n}. \quad (4.18)$$

On the other hand, by the definition (4.11) of $d_n$,
$$\Pr\left\{\frac{1}{n}\log\frac{1}{P_{V^n}(V^n)} \geq d_n + 2\gamma_n'\right\} \leq \delta_n. \tag{4.19}$$

Set $c_n = d_n + 2\gamma_n'$ and note that $\alpha_n \to 0$, $\delta_n \to 0$ $(n \to \infty)$ and $\gamma_n' = \dfrac{\gamma_n}{4}$, then by (4.18), (4.19) we have
$$\lim_{n\to\infty}\left(\Pr\left\{\frac{1}{n}\log\frac{1}{P_{V^n}(V^n)} \geq c_n\right\} \right.$$
$$\left. + \Pr\left\{\frac{1}{n}\log\frac{W^n(\tilde{Y}^n|\tilde{X}^n)}{P_{\tilde{Y}^n}(\tilde{Y}^n)} \leq c_n + \frac{1}{4}\gamma_n\right\}\right) = 0.$$

Finally, resetting $\tilde{X}^n\tilde{Y}^n$, $\dfrac{1}{4}\gamma_n$ as $X^nY^n$ and $\gamma_n$, respectively, we conclude that condition 2), i.e., (4.9) holds. □

Next, we show the equivalent information-spectrum form of the necessary condition (4.3) in Theorem 4.2.

**Theorem 4.4 (Equivalence of necessary conditions)** The following two conditions are equivalent:

1) For some channel input $\mathbf{X} = \{X^n\}_{n=1}^\infty$, which may be arbitrarily correlated to the output of the source $\mathbf{V} = \{V^n\}_{n=1}^\infty$, and for any sequence $\{\gamma_n\}_{n=1}^\infty$ satisfying condition (4.1), it holds that
$$\lim_{n\to\infty}\Pr\left\{\frac{1}{n}\log\frac{W^n(Y^n|X^n)}{P_{Y^n}(Y^n)} \leq \frac{1}{n}\log\frac{1}{P_{V^n}(V^n)} - \gamma_n\right\} = 0, \tag{4.20}$$

where $Y^n$ is the channel output via $W^n$ due to the channel input $X^n$.

2) (**Domination**) For some channel input $\mathbf{X} = \{X^n\}_{n=1}^\infty$, which may be arbitrarily correlated to the output of the source $\mathbf{V} = \{V^n\}_{n=1}^\infty$, and for some sequence $\{c_n\}_{n=1}^\infty$ and any sequence $\{\gamma_n\}_{n=1}^\infty$ satisfying condition (4.1), it holds that
$$\lim_{n\to\infty}\left(\Pr\left\{\frac{1}{n}\log\frac{1}{P_{V^n}(V^n)} \geq c_n\right\}\right.$$
$$\left. + \Pr\left\{\frac{1}{n}\log\frac{W^n(Y^n|X^n)}{P_{Y^n}(Y^n)} \leq c_n - \gamma_n\right\}\right) = 0, \tag{4.21}$$

where $Y^n$ is the channel output via $W^n$ due to the channel input $X^n$.

*Proof:* This theorem can be proved in the entirely same manner as in the proof of Theorem 4.3 with $\gamma_n$ replaced by $-\gamma_n$. □

**Remark 4.3** (*Separation Theorem*) The necessary condition 2) in Theorem 4.4 means that the entropy spectrum of tha source and the mutual information spectrum of the channel are asymptotically completely split except for the part of asymptotically vanishing width $\gamma_n$, and the former is placed to the left of the latter. This observation corresponds to the implication of the sufficient condition 2) in Theorem 4.3 (cf. Remark 4.2). If we disregard the asymptotically vanishing terms $\pm \gamma_n$, condition 2) in Theorem 4.3 together with condition 2) in Theorem 4.4 can be regarded as providing with a "necessary and sufficient condition" for the source $\mathbf{V} = \{V^n\}_{n=1}^{\infty}$ to be transmissible over the channel $\mathbf{W} = \{W^n\}_{n=1}^{\infty}$. Thus, in view of Remark 4.2, we can say that *Separation Theorem* continues to hold in a wider sense also for the general joint source-channel coding in consideration. □

**Remark 4.4** Actually, the definition of *domination* given by Vembu, Verdú and Steinberg [5] is not condition 2) in Theorem 4.4 but the following:

2′) (Domination) For *any* sequence $\{c_n\}_{n=1}^{\infty}$ and for any sequence $\{\gamma_n\}_{n=1}^{\infty}$ satisfying condition (4.1), there exists some channel input $\mathbf{X} = \{X^n\}_{n=1}^{\infty}$ such that
$$\lim_{n \to \infty} \left( \Pr\left\{\frac{1}{n} \log \frac{1}{P_{V^n}(V^n)} \geq c_n\right\} \right.$$
$$\left. \times \Pr\left\{\frac{1}{n} \log \frac{W^n(Y^n|X^n)}{P_{Y^n}(Y^n)} \leq c_n - \gamma_n\right\} \right) = 0 \quad (4.22)$$

holds, where $Y^n$ is the channel output via $W^n$ due to the channel input $X^n$.

It is easy to see that this necessary condition 2′) is implied by the necessary condition 2) in Theorem 4.4. In fact, the latter is much *stronger* than the former as necessary conditions for the transmissibility. □

## 5  $\varepsilon$-Transmissibility Theorem

So far we have considered only the case where the error probability $\varepsilon_n$ satisfies the condition $\lim_{n \to \infty} \varepsilon_n = 0$. However, we can relax this condition as

follows:
$$\limsup_{n\to\infty} \varepsilon_n \le \varepsilon, \qquad (5.1)$$
where $\varepsilon$ is any constant such that $0 \le \varepsilon < 1$. (It is obvious that the special case with $\varepsilon = 0$ coincides with the case that we have considered so far.) We now say that the source $\mathbf{V}$ is $\varepsilon$-*transmissible* over the channel $\mathbf{W}$ when there exists an $(n, \varepsilon_n)$ code satisfying condition (5.1).

Then, the same arguments as in the previous sections with due slight modifications lead to the following two theorems in parallel with Theorem 4.1 and Theorem 4.2, respectively:

**Theorem 5.1 ($\varepsilon$-Direct theorem)** Let $\mathbf{V} = \{V^n\}_{n=1}^{\infty}$, $\mathbf{W} = \{W^n\}_{n=1}^{\infty}$ be a general source and a general channel, respectively. If there exist some channel input $\mathbf{X} = \{X^n\}_{n=1}^{\infty}$, which may be arbitrarily correlated to the output of the source $\mathbf{V} = \{V^n\}_{n=1}^{\infty}$, and also some sequence $\{\gamma_n\}_{n=1}^{\infty}$ of real nubers with
$$\gamma_n > 0, \ \gamma_n \to 0 \text{ and } n\gamma_n \to \infty \quad (n \to \infty) \qquad (5.2)$$
for which it holds that
$$\limsup_{n\to\infty} \Pr\left\{\frac{1}{n}\log\frac{W^n(Y^n|X^n)}{P_{Y^n}(Y^n)} \le \frac{1}{n}\log\frac{1}{P_{V^n}(V^n)} + \gamma_n\right\} \le \varepsilon, \qquad (5.3)$$
then the source $\mathbf{V} = \{V^n\}_{n=1}^{\infty}$ is $\varepsilon$-transmissible over the channel $\mathbf{W} = \{W^n\}_{n=1}^{\infty}$, where $Y^n$ is the channel output via $W^n$ due to the channel input $X^n$. □

**Theorem 5.2 ($\varepsilon$-Converse theorem)** Suppose that a general source $\mathbf{V} = \{V^n\}_{n=1}^{\infty}$ is $\varepsilon$-transmissible over a general channel $\mathbf{W} = \{W^n\}_{n=1}^{\infty}$. Then, for some channel input $\mathbf{X} = \{X^n\}_{n=1}^{\infty}$, which may be arbitrarily correlated to the output of the source $\mathbf{V} = \{V^n\}_{n=1}^{\infty}$, and for any sequence $\{\gamma_n\}_{n=1}^{\infty}$ satisfying condition (5.2), it holds that
$$\limsup_{n\to\infty} \Pr\left\{\frac{1}{n}\log\frac{W^n(Y^n|X^n)}{P_{Y^n}(Y^n)} \le \frac{1}{n}\log\frac{1}{P_{V^n}(V^n)} - \gamma_n\right\} \le \varepsilon, \qquad (5.4)$$
where $Y^n$ is the channel output via $W^n$ due to the channel input $X^n$. □

It should be noted here that such a sufficient condition (5.3) as well as such a necessary condition (5.4) for the $\varepsilon$-transmissibility cannot actually be derived in the way of generalizing the strict domination in (4.9) and the domination in (4.21). It should be noted also that, under the $\varepsilon$-transmissibility criterion, *Separation Theorem* does *not* necessarily hold.

# 6 Separation Theorem for a Class of Sources and Channels

In this section, we consider, as a special case of Theorem 4.1~Theorem 4.4, the case where either the source $\mathbf{V} = \{V^n\}_{n=1}^{\infty}$ or the channel $\mathbf{W} = \{W^n\}_{n=1}^{\infty}$ satisfies the *strong converse* property. In this special case, Theorem 4.1 and Theorem 4.2 can be written in much simpler forms.

In order to show this, we need some preparations. Let $R_f(\mathbf{V})$, $C(\mathbf{W})$ denote the infimum of all achievable fixed-length coding rates for the source $\mathbf{V} = \{V^n\}_{n=1}^{\infty}$ and the capacity for the channel $\mathbf{W} = \{W^n\}_{n=1}^{\infty}$, respectively. A general source $\mathbf{V} = \{V^n\}_{n=1}^{\infty}$ is said to satisfy the strong converse property if the probability $\varepsilon_n$ of decoding error for fixed-length source coding with any rate $R$ such that $R < R_f(\mathbf{V})$ necessarily approaches one as $n$ tends to $\infty$ (cf. Han [9]). Moreover, a general channel $\mathbf{W} = \{W^n\}_{n=1}^{\infty}$ is said to satisfy the strong converse property if the probability $\varepsilon_n$ of decoding error for channel coding with any rate $R$ such that $R > C(\mathbf{W})$ necessarily approaches one as $n$ tends to $\infty$ (cf. Verdú and Han [8]).

Define[§]

$$\overline{H}(\mathbf{V}) \equiv \text{p-}\limsup_{n \to \infty} \frac{1}{n} \log \frac{1}{P_{V^n}(V^n)},$$

$$\underline{I}(\mathbf{X}; \mathbf{Y}) \equiv \text{p-}\liminf_{n \to \infty} \frac{1}{n} \log \frac{W^n(Y^n|X^n)}{P_{Y^n}(Y^n)},$$

where $Y^n$ is the channel output via $W^n$ due to the channel input $X^n$ and we have put $\mathbf{X} = \{X^n\}_{n=1}^{\infty}$, $\mathbf{Y} = \{Y^n\}_{n=1}^{\infty}$. Then, we have

**Theorem 6.1** (Han and Verdú [6])

$$R_f(\mathbf{V}) = \overline{H}(\mathbf{V}). \tag{6.1}$$

**Theorem 6.2** (Verdú and Han [8])

$$C(\mathbf{W}) = \sup_{\mathbf{X}} \underline{I}(\mathbf{X}; \mathbf{Y}). \tag{6.2}$$

---

[§]For any sequence $\{Z_n\}_{n=1}^{\infty}$ of real-valued random variables, we define the *limit superior in probability* (cf. Han and Verdú [6]) of $\{Z_n\}_{n=1}^{\infty}$ by $\text{p-}\limsup_{n \to \infty} Z_n = \inf\{\beta|\lim_{n \to \infty} \Pr\{Z_n > \beta\} = 0\}$. Also, we define the *limit inferior in probability* (cf. Han and Verdú [6]) of $\{Z_n\}_{n=1}^{\infty}$ by $\text{p-}\liminf_{n \to \infty} Z_n = \sup\{\alpha|\lim_{n \to \infty} \Pr\{Z_n < \alpha\} = 0\}$.

**Theorem 6.3** (Han [9]) The necessary and sufficient condition for the source $\mathbf{V}$ to satisfy the strong converse property is

$$\overline{H}(\mathbf{V}) = \underline{H}(\mathbf{V}), \tag{6.3}$$

where

$$\underline{H}(\mathbf{V}) \equiv \text{p-}\liminf_{n \to \infty} \frac{1}{n} \log \frac{1}{P_{V^n}(V^n)}. \tag{6.4}$$

**Theorem 6.4** (Verdú and Han [8]) The necessary and sufficient condition for the channel $\mathbf{W}$ to satisfy the strong converse property is

$$\sup_{\mathbf{X}} \underline{I}(\mathbf{X}; \mathbf{Y}) = \sup_{\mathbf{X}} \overline{I}(\mathbf{X}; \mathbf{Y}), \tag{6.5}$$

where

$$\overline{I}(\mathbf{X}; \mathbf{Y}) \equiv \text{p-}\limsup_{n \to \infty} \frac{1}{n} \log \frac{W^n(Y^n|X^n)}{P_{Y^n}(Y^n)}. \tag{6.6}$$

With these preparations, we have the following separation theorem for the class of sources and channels as stated above.

**Theorem 6.5 (Separation theorem)** Let either the source $\mathbf{V} = \{V^n\}_{n=1}^{\infty}$ or the channel $\mathbf{W} = \{W^n\}_{n=1}^{\infty}$ satisfy the strong converse property. Then, the following two statements hold:

1) If $R_f(\mathbf{V}) < C(\mathbf{W})$, then the source $\mathbf{V}$ is transmissible over the channel $\mathbf{W}$. In this case, we can separate the source coding and the channel coding (*Separation of codings*).

2) If the source $\mathbf{V}$ is transmissible over the channel $\mathbf{W}$, then it must hold that $R_f(\mathbf{V}) \leq C(\mathbf{W})$.

*Proof:*
1): In the proof of this part, the assumption of the strong converse property is not necessary. Since $R_f(\mathbf{V}) = \overline{H}(\mathbf{V})$, $C(\mathbf{W}) = \sup_{\mathbf{X}} \underline{I}(\mathbf{X}; \mathbf{Y})$ by Theorem 6.1 and Theorem 6.2, the inequality $R_f(\mathbf{V}) < C(\mathbf{W})$ implies that condition 2) in Theorem 4.3 holds for $\mathbf{X} = \{X^n\}_{n=1}^{\infty}$ attaining the supremum $\sup_{\mathbf{X}} \underline{I}(\mathbf{X}; \mathbf{Y})$ if we put, for example, $c_n = \frac{1}{2}(R_f(\mathbf{V}) + C(\mathbf{W}))$. Therefore, the source $\mathbf{V}$ is transmissible over the channel $\mathbf{W}$.

2): If the source $\mathbf{V}$ is transmissible over the channel $\mathbf{W}$, then condition 2) in Theorem 4.4 holds, and hence we have

$$\lim_{n\to\infty} \Pr\left\{\frac{1}{n}\log\frac{1}{P_{V^n}(V^n)} \geq c_n\right\} = 0. \tag{6.7}$$

First, suppose that the source $\mathbf{V}$ satisfies the strong converse property. Then, it follows from (6.7), the definition of $\underline{H}(\mathbf{V})$ and Theorem 6.3 that

$$\overline{H}(\mathbf{V}) = \underline{H}(\mathbf{V}) \leq \liminf_{n\to\infty} c_n.$$

Moreover, by virtue of condition 2) in Theorem 4.4 we have

$$\lim_{n\to\infty} \Pr\left\{\frac{1}{n}\log\frac{W^n(Y^n|X^n)}{P_{Y^n}(Y^n)} \leq c_n - \gamma_n\right\} = 0. \tag{6.8}$$

Then, it follows from (6.8) and the definition of $\underline{I}(\mathbf{X};\mathbf{Y})$ that

$$\underline{I}(\mathbf{X};\mathbf{Y}) \geq \liminf_{n\to\infty} c_n,$$

where we have put $\mathbf{X} = \{X^n\}_{n=1}^\infty$, $\mathbf{Y} = \{Y^n\}_{n=1}^\infty$. Hence, it is concluded that

$$R_f(\mathbf{V}) = \overline{H}(\mathbf{V}) \leq \underline{I}(\mathbf{X};\mathbf{Y}) \leq \sup_{\mathbf{X}} \underline{I}(\mathbf{X};\mathbf{Y}) = C(\mathbf{W}).$$

Next, suppose that the channel $\mathbf{W}$ satisfies the strong converse property. Then, it follows from (6.8), the definition of $\overline{I}(\mathbf{X};\mathbf{Y})$ and Theorem 6.2, Theorem 6.4 that

$$\limsup_{n\to\infty} c_n \leq \overline{I}(\mathbf{X};\mathbf{Y}) \leq \sup_{\mathbf{X}} \overline{I}(\mathbf{X};\mathbf{Y}) = \sup_{\mathbf{X}} \underline{I}(\mathbf{X};\mathbf{Y}) = C(\mathbf{W}).$$

On the other hand, from (6.7) and the definition of $\overline{H}(\mathbf{V})$ we have

$$\overline{H}(\mathbf{V}) \leq \limsup_{n\to\infty} c_n.$$

Thus, it is again concluded that $R_f(\mathbf{V}) = \overline{H}(\mathbf{V}) \leq C(\mathbf{W})$. $\square$

**Remark 6.1** In the proof of Theorem 6.5 2) we have invoked the domination in the sense of Theorem 4.4 2). Notice, however, that the domination in the sense of Remark 4.4 2′) cannot lead us to establish Theorem 6.5 2).

**Example 6.1** If neither the source $\mathbf{V} = \{V^n\}_{n=1}^{\infty}$ nor the channel $\mathbf{W} = \{W^n\}_{n=1}^{\infty}$ satisfies the strong converse property, the statement 2) in Theorem 6.5 does *not* necessarily hold. For example, let $\mathcal{V} = \mathcal{X} = \mathcal{Y} = \{0,1\}$, and consider the case where, if $n$ is even, then $P_{V^n}$ is the uniform distribution on $\mathcal{V}^n$ and $W^n$ is the identity mapping; otherwise if $n$ is odd, then $P_{V^n}(0^n) = 1$, $W^n(0^n|\mathbf{x}) = 1$ ($\forall \mathbf{x} \in \mathcal{X}^n$), where $0^n$ denotes the sequence of $n$ 0's. It is obvious in this case that the source $\mathbf{V}$ is transmissible over the channel $\mathbf{W}$ with *zero* error probability, although neither the source $\mathbf{V}$ nor the channel $\mathbf{W}$ satisfies the strong converse property. However, since it is easy to check that $R_f(\mathbf{V}) = \log 2$ and $C(\mathbf{W}) = 0$, the statement 2) in Theorem 6.5 does not hold for this example. □

Let us now compare Theorem 6.5 with the standard classical separation theorem. To this end, we need to record some definitions as follows: A general source $\mathbf{V} = \{V^n\}_{n=1}^{\infty}$ is said to be *information-stable* (cf. Dobrushin [2], Pinsker [3]) if

$$\frac{\frac{1}{n}\log\frac{1}{P_{V^n}(V^n)}}{H_n(V^n)} \to 1 \quad \text{in prob.,} \tag{6.9}$$

where $H_n(V^n) = \frac{1}{n}H(V^n)$ and $H(V^n)$ stands for the entropy of $V^n$ (cf. Cover and Thomas [13]). Moreover, a general channel $\mathbf{W} = \{W^n\}_{n=1}^{\infty}$ is said to be *information-stable* (cf. Dobrushin [2], Hu [4]) if there exists a channel input $\mathbf{X} = \{X^n\}_{n=1}^{\infty}$ such that

$$\frac{\frac{1}{n}\log\frac{W(Y^n|X^n)}{P_{Y^n}(Y^n)}}{C_n(W^n)} \to 1 \quad \text{in prob.,} \tag{6.10}$$

where

$$C_n(W^n) = \sup_{X^n} \frac{1}{n} I(X^n; Y^n),$$

and $Y^n$ is the channel output via $W^n$ due to the channel input $X^n$; and $I(X^n; Y^n)$ is the mutual information between $X^n$ and $Y^n$ (cf. Cover and Thomas [13]). Then, we can summarize the typical classical separation theorem as follows, which is slightly different from Theorem 6.5:

**Theorem 6.6** (Vembu, Verdú and Steinberg [5]) Let the channel $\mathbf{W} = \{W^n\}_{n=1}^{\infty}$ be information-stable and suppose that the limit $\lim_{n\to\infty} C_n(W^n)$ exists. Or, let the source $\mathbf{V} = \{V^n\}_{n=1}^{\infty}$ be information-stable and suppose that the limit $\lim_{n\to\infty} H_n(V^n)$ exists. Then, the following two statements hold:

1) If $R_f(\mathbf{V}) < C(\mathbf{W})$, then the source $\mathbf{V}$ is transmissible over the channel $\mathbf{W}$. In this case, we can separate the source coding and the channel coding (*Separation of codings*).

2) If the source $\mathbf{V}$ is transmissible over the channel $\mathbf{W}$, then it must hold that $R_f(\mathbf{V}) \leq C(\mathbf{W})$.

It is not difficult to verify that, in the case where both of channel input alphabet $\mathcal{X}$ and channel output alphabet $\mathcal{Y}$ are non-finite *abstract* sets, either the strong converse property of the channel $\mathbf{W}$ assumed in Theorem 6.5 or the information-stability of the channel $\mathbf{W}$ together with the existence of the limit assumed in Theorem 6.6 is *not* implied by the other. In this sense, both theorems have their own rights. On the other hand, with *countably infinite* source alphabet $\mathcal{V}$, the information-stability of the source $\mathbf{V}$ together with the existence of the limit assumed in Theorem 6.6 implies the strong converse property of the source $\mathbf{V}$ assumed in Theorem 6.5. In this sense, Theorem 6.5 is *stronger* than Theorem 6.6.

It should be pointed out here that in general it is more or less easier to check the validity of the condition in Theorem 6.5 than that of the condition in Theorem 6.6.

Finally, let us now consider to generalize both of Theorem 6.5 and Theorem 6.6. In fact, we can strengthen Theorem 6.5 so as to include also Theorem 6.6 as a special case. To do so, let us define the concept of *semi-strong converse* property as follows. A general source $\mathbf{V} = \{V^n\}_{n=1}^{\infty}$ is said to satisfy the *semi*-strong converse property if for all divergent subsequences $\{n_i\}_{n=1}^{\infty}$ of positive integers such that $n_1 < n_2 < \cdots \to \infty$ it holds that

$$\text{p-}\limsup_{i \to \infty} \frac{1}{n_i} \log \frac{1}{P_{V^{n_i}}(V^{n_i})} = \overline{H}(\mathbf{V}). \tag{6.11}$$

A general channel $\mathbf{W} = \{W^n\}_{n=1}^{\infty}$ is said to satisfy the *semi*-strong converse property if for any channel input $\mathbf{X} = \{X^n\}_{n=1}^{\infty}$ and for all divergent subsequences $\{n_i\}_{n=1}^{\infty}$ of positive integers such that $n_1 < n_2 < \cdots \to \infty$ it holds that

$$\text{p-}\liminf_{i \to \infty} \frac{1}{n_i} \log \frac{W^{n_i}(Y^{n_i}|X^{n_i})}{P_{Y^{n_i}}(Y^{n_i})} \leq \sup_{\mathbf{X}} \underline{I}(\mathbf{X};\mathbf{Y}), \tag{6.12}$$

where $Y^n$ is the channel output via $W^n$ due to the channel input $X^n$.

With these definitions, we have the following separation theorem.

**Theorem 6.7 (Separation theorem)** Let either the source $\mathbf{V} = \{V^n\}_{n=1}^{\infty}$ or the channel $\mathbf{W} = \{W^n\}_{n=1}^{\infty}$ satisfy the *semi*-strong converse property. Then, the following two statements hold:

1) If $R_f(\mathbf{V}) < C(\mathbf{W})$, then the source $\mathbf{V}$ is transmissible over the channel $\mathbf{W}$. In this case, we can separate the source coding and the channel coding (*Separation of codings*).

2) If the source $\mathbf{V}$ is transmissible over the channel $\mathbf{W}$, then it must hold that $R_f(\mathbf{V}) \leq C(\mathbf{W})$.

*Proof:*

The proof is the same as that of Theorem 6.5 based on Theorem 4.4 2), except for that, here, owing to the assumed conditions, we directly have $\overline{H}(\mathbf{V}) \leq \liminf_{n \to \infty} c_n$ (if the source $\mathbf{V}$ satisfies the semi-strong converse property) and $\limsup_{n \to \infty} c_n \leq \sup_{\mathbf{X}} \underline{I}(\mathbf{X};\mathbf{Y})$ (if the channel $\mathbf{W}$ satisfies the semi-strong converse property). $\square$

**Remark 6.2** It is obvious that the strong converse property implies the *semi*-strong converse property for either the source or the channel. Therefore, Theorem 6.7 includes Theorem 6.5 as a special case. Similarly, it is not difficult to check that the information stability together with the existence of the limit implies the *semi*-strong converse property for either the source or the channel. Hence, Theorem 6.7 includes Theorem 6.6 as a special case. Thus, Theorem 6.7 is the *strongest* among these three separation theorems in the traditional sense.

Csiszár and Körner [12] have posed two operational standponits in source coding and channel coding, i.e., the *pessimistic standpint* and the *optiimistic standpint*. In their terminology, for source coding, the semi-strong convserse property is equivalent to the condition that both of the pessimistic standpint and the optiimistic standpint result in the same infimum of all achievable fixed-length source coding rates. Similarly, for channel coding, the semi-strong convserse property is equivalent to the condition that both of the pessimistic standpint and the optiimistic standpint result in the same supremum of all achievable channel coding rates. $\square$

**Example 6.2** Let us consider two different stationary memoryless sources $\mathbf{V}_1 = \{V_1^n\}_{n=1}^{\infty}$, $\mathbf{V}_2 = \{V_2^n\}_{n=1}^{\infty}$ with countably infinite source alphabet $\mathcal{V}$, and define its *mixed* source $\mathbf{V} = \{V^n\}_{n=1}^{\infty}$ by

$$P_{V^n}(\mathbf{v}) = \alpha_1 P_{V_1^n}(\mathbf{v}) + \alpha_2 P_{V_2^n}(\mathbf{v}) \quad (\mathbf{v} \in \mathcal{V}^n),$$

where $\alpha_1, \alpha_2$ are positive constants such that $\alpha_1 + \alpha_2 = 1$. Then, this mixed source $\mathbf{V} = \{V^n\}_{n=1}^{\infty}$ satisfies the semi-strong converse property but neither the strong converse property nor the information-stability.

Similarly, let us consider two different stationary memoryless channels $\mathbf{W}_1 = \{W_1^n\}_{n=1}^{\infty}$, $\mathbf{W}_2 = \{W_2^n\}_{n=1}^{\infty}$ with arbitrary abstract input and output alphabets $\mathcal{X}, \mathcal{Y}$, and define its *mixed* channel $\mathbf{W} = \{W^n\}_{n=1}^{\infty}$ by

$$W^n(\mathbf{y}|\mathbf{x}) = \alpha_1 W_1^n(\mathbf{y}|\mathbf{x}) + \alpha_2 W_2^n(\mathbf{y}|\mathbf{x}) \quad (\mathbf{x} \in \mathcal{X}^n, \mathbf{y} \in \mathcal{Y}^n).$$

Then, this mixed channel $\mathbf{W} = \{W^n\}_{n=1}^{\infty}$ satisfies the semi-strong converse property but neither the strong converse property nor the information-stability. $\square$

# References

[1] C. E. Shannon, "A mathematical theory of communication," *Bell System Technical Journal*, vol.27, pp.379-423, pp. 623-656, 1948

[2] R. L. Dobrushin, "A general formulation of the fundamental Shannon theorem in information theory," *Uspehi Mat. Acad. Nauk. SSSR*, vol.40, pp.3-104, 1959: Translation in *Transactions of American Mathematical Society*, Series 2, vol.33, pp.323-438, 1963

[3] M. S. Pinsker, *Information and Information Stability of Random Variables and Processes*, Holden-Day, San Francisco, 1964

[4] G. D. Hu, "On Shannon theorem and its converse for sequence of communication schemes in the case of abstract random variables," in *Trans. 3rd Prague Conference on Information Theory, Statistical Decision Functions, Random Processes*, Czechslovak Academy of Sciences, Prague, pp. 285-333, 1964

[5] S. Vembu, S. Verdú and Y. Steinberg, "The source-channel separation theorem revisited," *IEEE Transactions on Information Theory*, vol.IT-41, no.1, pp. 44-54, 1995

[6] T.S. Han and S. Verdú, "Approximation theory of output statistics," *IEEE Transactions on Information Theory*, vol.IT-39, no.3, pp. 752-772, 1993

[7] S. Verdú and T. S. Han, "The role of the asymptotic equipartition property in noiseless source coding," *IEEE Transactions on Information Theory*, vol.IT-43, no.3, pp.847-857, 1997

[8] S. Verdú and T.S. Han, "A general formula for channel capacity," *IEEE Transactions on Information Theory*, vol.IT-40, no.4, pp.1147-1157, 1994

[9] T. S. Han, *Information-Spectrum Methods in Information Theory*, Baifukan-Press, Tokyo, 1998 (in Japanese)

[10] T. S. Han, "The reliability functions of the general source with fixed-length coding," to appear in *IEEE Transactions on Information Theory*

[11] A. Feinstein, "A new basic theorem of information theory," *IRE Trans. PGIT*, vol.4, pp.2-22, 1954

[12] I. Csiszár and J. Körner, *Information Theory: Coding Theorems for Discrete Memoryless Systems*, Academic Press, New York, 1981

[13] T. M. Cover and J. A. Thomas, *Elements of Information Theory*, Wiley, New York, 1991